\documentclass[a4paper,12pt,leqno]{article} 

\usepackage{lineno}

\usepackage{amsmath,amsfonts} 
\usepackage{graphics}
\usepackage[T1]{fontenc}
\usepackage{lineno}
\usepackage[utf8]{inputenc}
\usepackage{amsmath,mathrsfs,amssymb,amsfonts,nicefrac}
\usepackage[english]{babel}
\usepackage[active]{srcltx}
\usepackage{tikz}
\usepackage{dsfont}
\usetikzlibrary{patterns}

\numberwithin{equation}{section}

\newtheorem{thm}{Theorem}[section]   
\newtheorem{cor}[thm]{Corollary}   
\newtheorem{prop}[thm]{Proposition}   
\newtheorem{lm}[thm]{Lemma}  

\newtheorem{rem}[thm]{Remark}  

\newcommand{\RR}{\mathbb{R}}   
\newcommand{\R}{\mathbb{R}}   
\newcommand{\di}{\displaystyle}

\renewcommand{\epsilon}{\varepsilon}
\newcommand{\e}{\varepsilon}
\def\un{{\mathbf{1}}}

\begin{document}   
   
\title{\textbf{The leading edge of a free boundary interacting  with a line of fast diffusion}}

\author{
{\bf Luis A. Caffarelli} \\ 
The University of Texas at Austin\\
Mathematics Department RLM 8.100, 2515 Speedway Stop C1200\\
 Austin, Texas 78712-1202, U.S.A.  \\
\texttt{caffarel@math.utexas.edu}
\\[2mm]
{\bf Jean-Michel Roquejoffre} \\ 
Institut de Math\'ematiques de Toulouse (UMR CNRS 5219) \\ 
Universit\'e Toulouse III,
118 route de Narbonne\\
31062 Toulouse cedex, France \\ 
\texttt{jean-michel.roquejoffre@math.univ-toulouse.fr}}
  \date{}
\maketitle   
\begin{abstract} \noindent The goal of this work is to explain an unexpected feature of the expanding level sets of the solutions of a system where a half plane, in which reaction-diffusion phenomena take place, exchanges mass with a line having a large diffusion of its own. The system was proposed by H. Berestycki, L. Rossi and the second author \cite{BRR2} as a model of enhancement of biological invasions by a line of fast diffusion. It was observed numerically by A.-C. Coulon \cite{ACC} that the leading edge of the front, rather than being located on the line, was in the lower half plane. 

\noindent We explain this behaviour for a free boundary problem much related to the system for which the simulations were made. We construct travelling waves for this problem, analyse their free boundary near the line, and prove that it has the behaviour predicted by the numerical simulations.
\end{abstract}

\section{Introduction}    

\subsection{Model and question}
 Consider the cylinder
$
\Sigma=\{(x,y)\in\RR\times(-L,0)\}.
$
We look for a real $c>0$, a function $u(x)$, defined for $x\in\RR$, a function $v(x,y)$ defined in $\Sigma$, and a curve $\Gamma\subset\Sigma$ such that
\begin{equation}
\label{e1.200}
\left\{
\begin{array}{rll}
-d\Delta v+c\partial_x v=&0\quad(x,y)\in\{v>0\}\\
\vert\nabla v\vert=&1\quad((x,y)\in\Gamma:=\Sigma\cap\partial\{v>0\}\\
\ \\
-Du_{xx}+c\partial_xu+1/\mu u-v=&0\quad\hbox{for $x\in\RR$, $y=0$}\\
v_y=&\mu u-v\quad\hbox{for $x\in\RR$, $y=0$ and $v(x,0)>0$}\\
u_y(x,-L)=&0\\
u(-\infty)=1/\mu,\ u(+\infty)=&0,\quad v(-\infty,y)=1,\ v(+\infty,y)=0.
\end{array}
\right.
\end{equation}
In \eqref{e1.200}, the real numbers $\mu,d,D$ are fixed positive constants, and the problem inside $\Sigma$ is a well-known free boundary problem. We will also consider then the following more compact problem, with unknowns $(c,\Gamma, u)$, the function $u$ being this time defined in $\overline\Sigma$, solving
\begin{equation}
\label{e1.2}
\left\{
\begin{array}{rll}
-d\Delta u+c\partial_x u=&0\quad(x,y)\in\{u>0\}\\
\vert\nabla u\vert=&1\quad((x,y)\in\Gamma:=\Sigma\cap\partial\{u>0\}\\
\ \\
-Du_{xx}+c\partial_xu+1/\mu u_y=&0\quad\hbox{for $x\in\RR$, $y=0$}\\
u_y(x,-L)=&0\\
u(-\infty,y)=&1,\quad u(+\infty,y)=0.
\end{array}
\right.
\end{equation}
In both problems \eqref{e1.200} and \eqref{e1.2}, we will see that $v_x\leq 0$ (resp. $u_x\leq 0$) inside $\Sigma$, 
and that the free boundary $\Gamma$ inside $\Sigma$ will be an analytic curve. Assume that $\Gamma$ intersects the line $\{y=0\}$, say at $(x,y)=(0,0)$. We ask for the behaviour of $\varphi$ near $y=0$. 

\subsection{Motivation} 
Our starting point is the system proposed by H. Berestycki, L. Rossi and the second author to model the speed-up of biological invasions by lines of fast diffusion in \cite{BRR2}. In this model, the two-dimensional lower half-plane ("`the field"'), in which reaction-diffusion phenomena occur, interacts with the $x$ axis ("`the road") which has a much faster diffusion $D$ of its own. It will sometimes be useful to assume $D\geq d$, but not always. Call $u(t,x)$ the density of individuals on the road, and $v(t,x,y)$ the density of individuals in the field. The road yields the fraction $\mu u$ to the field, and retrieves the fraction $\nu v$ in exchange; the converse process occurs for the field.  The system for $u$ and $v$ is
\begin{equation}
\label{e1.3}
\begin{array}{rll}
\partial_t u-D \partial_{xx} u= & \nu v(t,x,0)-\mu u \   \    &x\in\R\\
\partial_t v-\Delta v=&f(v)\    \    &(x,y)\in\R\times\R_-\\
\partial_y v(t,x,0)=&\mu u(t,x,t)-  \nu v(t,x,0)\   \    &x\in\R.
\end{array}
\end{equation}
Here $f$ is the usual logistic term, $f(v)=v-v^2$. A model involving only the unknown $u$ can be obtained by forcing the (biologically reasonable) equality
$\phi(x)=\psi(x,0)$; in other words we do the (formal)  limit $\delta\to0$ of $\nu=\mu=\di\frac1\delta$. Still arguing in a formal way, we obtain $u=v$ on the road, and the exchange term is simply $v_y$. Thus, we obtain
\begin{equation}
\label{e1.3000}
\begin{array}{rll}
\partial_t v-D \partial_{xx} v+v_y(t,x,0)= & 0 \   \    &x\in\R\\
\partial_t v-\Delta v=&f(v)\    \    &(x,y)\in\R\times\R_-\\
\end{array}
\end{equation}
\noindent From now on, as is rather intuitively clear from the biological modelling, we will call \eqref{e1.200} the two-species model, whereas problem \eqref{e1.2} will be tge single-species model.
For the time being, let us only argue on system \eqref{e1.3}. The first question is how the stable state $(\nu/\mu,1)$ invades the unstable state $(0,0)$. In \cite{BRR4} it is computed with $o_{t\to+\infty}(1)$ precision: for each direction $e$ in the field, the level sets of $v$ move with a velocity $w_*(e)$ which, quite surprisingly,  does not obey the Huygens principle. The next step is to describe the asymptotic level sets with $O_{t\to+\infty}(1)$ precision; for this purpose numerical simulations were carried out by A.-C. Coulon, the simulations, which are part of a larger program of her thesis \cite{ACC}. 
The above figures account for some of her results; the parameters are 
$$f(v)=v-v^2,\ \ D=10,\ \ u(0,x)={\bf 1}_{[-1,1]}(x),\ \ v(0,x,y)\equiv0.
$$
 The top figure represents the levels set 0.5 of $v$ at times 10, 20, 30, 40; the bottom figure represents the shape of $v(40,x,y)$.

\includegraphics{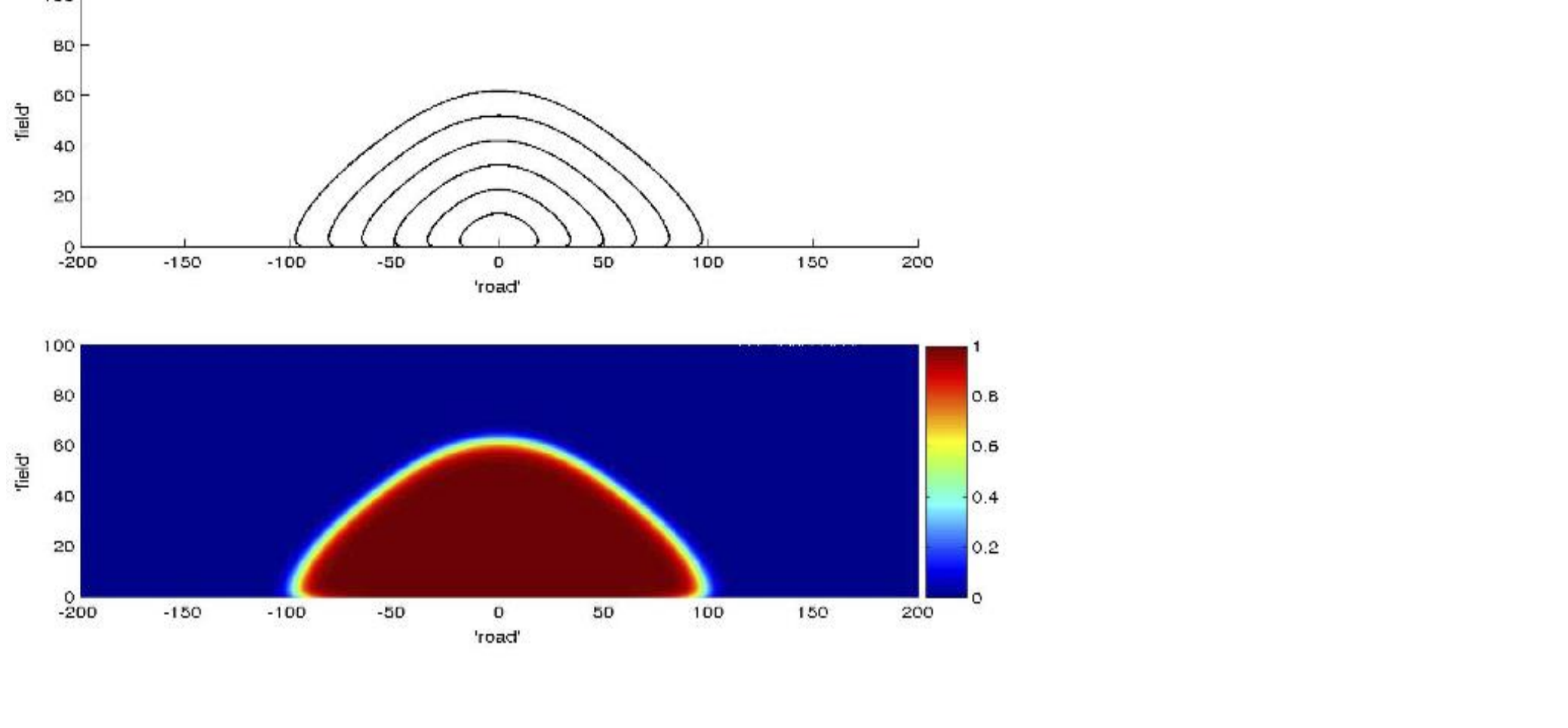}		

\noindent We have found these figures surprising, all the more as they are quite robust with respect to all the parameters. Indeed, a naive intuition would suggest that the leading edge of the invasion is located on the road, especially for large $D$. Such is manifestly not the case, the leading edge appears to be located in the field, at a distance to the road which seems to remain more or less constant in time. A heuristic explanation is the following: the term  $v(t,x,0)-\mu u(t,x)$ acts as an effective reaction term for $u$; given that everything suggests that the invasion is driven by the road, especially for the large values of $D$.  The immediate consequence is that $v_y(t,x,0)<0$, so the function $v$ increases in the vicinity of the road, hence the observed behaviour. The goal of this paper is to find a mathematically rigorous account of it. 

\noindent Working directly on \eqref{e1.3} or \eqref{e1.3000} to explain the above simulations appeared to be a difficult task. To circumvent the difficulty, we choose here  to work on a limiting model, which will keep the main features of \eqref{e1.3} or \eqref{e1.3000}, this is why we have come up with \eqref{e1.200} and \eqref{e1.2}. This system is indeed a limiting model of \eqref{e1.3}, let us explain why. Before that, we consider a permanent regime $(\partial_tu,\partial_tv)=(0,0)$ of \eqref{e1.3}, it is logical to think that it will take the form of a travelling wave 
\begin{equation}
\label{e2.30}
(u(t,x),v(t,x,y))=(\phi(x+ct),\psi(x+ct,y)),\quad c>0.
\end{equation}
It is proved in \cite{DiR} that travelling waves are indeed attracting for \eqref{e1.2} (and, more interestingly, the paper describes how the convergence occurs when $D$ is very large). We have  consider a propagation from right to left, thus the couple   $(\phi,\psi)$ is an orbit of the following system
\begin{equation}
\label{e1.4}
\begin{array}{rll}
-D \partial_{xx} \phi+c\phi_x= & \nu\psi(x,0)-\mu \phi \   \    &x\in\R\\
-\Delta \psi+c\psi_x=&f(\psi)\    \    &(x,y)\in\Sigma\\
\partial_y \psi(x,0)=&\mu \phi(x)-  \nu\psi(x,0)\   \    &x\in\R\\
\end{array}
\end{equation}
where we have omitted the Neumann condition for short. As for model \eqref{e1.3000}, the corresponding system is (we look for $v(t,x,y)$ under the form $\psi(x+ct,y)$):
\begin{equation}
\label{e1.6}
\begin{array}{rll}
-D \partial_{xx} \psi+c\phi_x+\psi_y= & 0 \   \    &x\in\R\\
-\Delta \psi+c\psi_x=&f(\psi)\    \    &(x,y)\in\Sigma\\
\end{array}
\end{equation}
A free boundary problem is obtained (once again in a formal way) in the limit $\e\to0$ of a sequence of solutions $(u_\e,v_\e)$ to \eqref{e1.4} with  $f=f_\e$, an approximation of the Dirac mass  $\delta_{\psi=1}$; if $\Gamma$ is the unknown boundary $\partial(\{\psi=1\})$, and if we set $u=1/\mu-\phi$, $v:=1-\psi$, we obtain \eqref{e1.200}. To retrieve \eqref{e1.2}, it suffices to set $u:=1-\psi$ in \eqref{e1.6}.

\noindent  The question is therefore whether, and how, the free boundary $\Gamma$ meets the fixed boundary $\{y=0\}$. For the two unknowns problem \eqref{e1.200}, the scaling - both inside the domain and at the vicinity of the hitting point - is Lipschitz, which allows the use of a large body of existing ideas. On the other hand, Model \eqref{e1.2} is not of the standard type, because the characteristic scales around the free boundary are different inside $\Sigma$ and on the fixed boundary are different: in the latter it is $y\sim x^2$.

\subsection{Results}
Before studying free boundary $\Gamma$ of Problem \eqref{e1.2} near the origin, we first should make sure that it has a solution. Cleary, the issue is what happens in the vicinity of the axis $\{y=0\}$, therefore we start from a situation for which travelling waves are known to exist. 

\begin{thm}
\label{t2.1}
System \eqref{e1.2} has a solution $(c,\Gamma,u,v)$. We have $c>0$ and $\partial_xu\leq 0$ (and $<0$ to the left of $\Gamma$). The function $v$ is globally Lipschitz: $\vert\nabla v\vert\leq C$ for some universal $C$ (therefore $u$ is $C^{1,1}$ on the line). The free boundary $\Gamma$ is a graph $(\varphi(y),y)$, and also an analytic curve in $\{y<0\}$. Moreover it intersects the $x$-axis, so we may choose $\varphi(0)=0$. 
\end{thm}
\noindent {\bf Remark 1.} {\it Uniqueness (up to translations) of $(c,\Gamma,u,v)$ is probably true. The only issue is to examine the behaviour of two solutions whose contact point lies at the intersection of the $x$ axis and their respective free boundaries. This will not be pursued here.}

\noindent Next, we study the free boundary in a neighbourhood of the origin.
\begin{thm}
\label{t2.2}
The free boundary $\Gamma$ hits the line $\{y=0\}$ at a point where $u>0$. We may assume that the hitting point is $(0,0)$. Moreover, in a neighbourhood of $(0,0)$ $\Gamma$ is a graph in the $x$ variable $y=\psi(x),\ x<0,$, and there is $\gamma>0$ such that: 
\begin{equation}
\label{e2.10001}
\psi(x)=\gamma x+o_{x\to0^+}(x).
\end{equation}
\end{thm}

\noindent Let us turn to the model with one species.
\begin{thm}
\label{t2.3}
Assume $D\geq d$. System \eqref{e1.2} has a solution $(c,\Gamma,u)$. We have $c>0$ and $\partial_xu\leq 0$ (and $<0$ to the left of $\Gamma$). The function $u$ is globally Lipschitz: $\vert\nabla u\vert\leq C$ for some universal $C$. The free boundary $\Gamma$ is a graph $(\varphi(y),y)$, and also an analytic curve in $\{y<0\}$. Moreover it intersects the $x$-axis, so we may choose $\varphi(0)=0$. 
\end{thm}
\noindent {\bf Remark 2.} {\it We do not know whether the assumption $D\geq d$ is indispensable, or just technical. In any case, it is consistent with our goal to study phenomena driven by a large diffusion on the road. In \cite{BRR2}, the significant threshold is $D=2d$, this is where the velocity of the wave for \eqref{e1.2} exceeds that of the plane Fisher-KPP wave.}

\begin{thm}
\label{t2.4}
Assume $D\geq d$. In a neighbourhood of $(0,0)$ $\Gamma$ is a graph in the $x$ variable: 
\begin{equation}
\label{e2.1000}
y=\psi(x),\ x<0,\quad\quad\hbox{with}\ \psi(x)=-\frac{x^2}{2D}+o_{x\to0^+}(x^2).
\end{equation}
\end{thm}
\subsection{Discussion, organisation of the paper}
We first notice, as far as the two species is concerned, an interesting loss of boundary condition between the reaction-diffusion system \eqref{e1.4} and the free boundary problem \eqref{e1.2}, to the right of the right of the free boundary. In other words, in this area, the road does not exchange individuals anymore with the field. Thus there is an asymptotic decoupling between the value of $u$ and the value of $v$ at the free boundary, and accounts quite well for the numerical simulations. Let us revert to the old unknowns $\phi$ and $\psi$, which denote respectively the density of individuals on the road and in the field. At the intersection between the invasion front and the road, the density of individuals is at its maximum, whereas the road keeps feeding the field with individuals, as is stated by the exchange condition 
$$
\psi_y\sim \mu \phi-\psi\sim \mu u(0,0)>0.
$$
Therefore, the invasion front in the field can only go further, which explains that its leading edge is not located on the road. This explains the simulations in \cite{ACC}.

\noindent Next, we observe a different behaviour for the one species model. A heuristic reason is that the model is trying to accomodate both the exchange condition and the free boundary condition, and this is only done at the expense of a breakdown of the homogeneity near the road. The resulting sitation is the interaction of an obstacle problem on the road and a solution of the one phase problem in the field. This observation will be investigated further in a subsequent work.

\noindent The last remark concerns the derivation of the one species model from the two species model, in the limit of infinite exchange terms $\mu=\nu=\di\frac1\delta$. This passage is done in a rigorous way in \cite{Di1} in the framework of a reaction-diffusion system; it would probably be possible to do it, such a task would probably be not entirely trivial. We postpone this matter to a future work, as it would not really add more to the understanding of the question at stake here.

\medskip
\noindent The  paper is organised as follows. In Section 2, we construct a solution to system \eqref{e1.200}. This is done by a classical approximation by a family of semilinear equations where the nonlinearity$f_\e(u)$ converges, in the measure sense, to $\delta_{\partial_{u>0}}$. The idea, as well as the inspiration for the proof, is taken from \cite{BCN}. In Section 3, we study the free boundary of the two species model, the main argument will be a Liouville type theorem for a special class of global solutions. In Section 4, we construct the travelling wave, and the main part of the analysis is the gradient bound for $u$. Finally, in Section 5, we prove Theorem \ref{t2.4}. The chief argument in all these (at times technical) considerations is  that the free boundary condition creates such a rigidity that, unless the solutions behave as they are expected to behave, basic properties such as, for instance, the positivity of $u$, will not hold true.  
\section{The travelling wave in the two species model}
Let us consider a smooth function $\varphi(u)$, defined on $\RR_+$, positive on $[0,1)$, zero outside, and such that
\begin{equation}
\label{e3.1117}
\int_0^{+\infty}u\varphi(u)du=\frac1{2d}.
\end{equation}
Consider the sequence of reaction terms
\begin{equation}
\label{e3.1}
f_\e(u)=\frac{u}{\e^2}\varphi(\frac{u}\e),
\end{equation}
We will obtain a solution to \eqref{e1.200} as the limit, as $\e\to0$, of a sequence $(c_e,u_\e, v_\e)_\e$ of solutions of 
\begin{equation}
\label{e3.2000}
\left\{
\begin{array}{rll}
-d\Delta v+c\partial_x v+f_\e(v)=&0\quad(x,y)\in\Sigma\\
\ \\
-Du_{xx}+c\partial_xu+1/\mu u-v=&0\quad\hbox{for $x\in\RR$, $y=0$}\\
v_y=&\mu u-v\quad\hbox{for $x\in\RR$, $y=0$}\\
v_y(x,-L)=&0\\
\ \\
v(-\infty,y)=1,\ v(+\infty,y)=&0,\quad u(-\infty)=1/\mu,\ u(+\infty)=0.
\end{array}
\right.
\end{equation}
This is by no means a new idea. Such reaction terms, which date back to Zeldovich \cite{ZBLM} have been quite useful to study asymptotic models in flame propagation, that we will not quote here because of their abundance. The first rigorous passage from the reaction-diffusion equation to the free boundary problem is done in the context of 1D travelling waves, by Berestycki, Nicolaenko and Scheurer \cite{BNS}. In several space dimensions, it was done by Berestycki, Nirenbrg and the first author \cite{BCN}, a work that we will much use here. Let us briefly recall why a limiting solution $(c,u,v)$ to \eqref{e3.2000} will develop a free boundary. For this, we consider the simple one-dimensional problem 
\begin{equation}
\label{e3.3}
\begin{array}{rll}
-du''+cu'+f_\e(u)=&0\ \hbox{on $\RR$})\\
 u(-\infty)=0,\ u(+\infty)&=1.
\end{array}
\end{equation}
Standard arguments show that, for a solution $(c_\e,u_\e)$ to \eqref{e3.3}, we have $c_\e>0$ and $u_\e'<0$. We may always assume that $u_\e(0)=\e$, so, for $x<0$ we have $u_\e(x)=1-(1-\e)e^{cx}$. As for $x>0$ we set
$$
\xi=\frac{x}\e,\quad u_\e(x)=\e p_\e(\frac{x}\e).
$$ 
The function $p_\e$ solves
$$
-dp''+\e c_\e p'+p\varphi(p)=0\ \hbox{on}\ \RR_+,\quad p(0)=1,\ p(+\infty)=0.
$$
Once again, standard arguments show that the term $\e c_\e p'$ may be neglected, so that multiplication by $p'$ and integration over $\RR_+$ yields
$$
d\frac{(p_\e')^2(0)}2=\int_0^{+\infty}p\varphi(p)dp=\frac1{2d}.
$$ 
Scaling back, we obtain $u_\e(0)\sim1$, matching derivatives yields $c=1$. The limit $(c,u)$ of $(u_\e,c_\e)$ satisfies therefore 
$$
-du''+cu'=0\ \hbox{on}\ \RR_-,\quad d[u'](0)=1.
$$
This is the one-dimensional version of the problem inside $\Sigma$; here, $\Gamma$ is the point $x=0$.

Let us come back to \eqref{e3.2000}. For every $\e>0$, \eqref{e3.2000} has (Dietrich \cite{Di1}, Theorem 1) a unique solution $(c_\e,u_\e)$ such that $c_\e>0$, $0<u_\e<1$ and $\partial_xu_\e,\partial_xv_\e<0$; we will show that, up to a subsequence, $(c_\e,u_\e)$ converges, as $\e\to0$, to a solution of \eqref{e1.2}. 
\subsection{Basic bounds}
The first task is to show that the travelling wave velocity is uniformly bounded from above, that is
\begin{prop}
\label{p3.1000}
There is $K>0$ independent of $\e$ such that $c_\e\leq K$.
\end{prop}
\noindent{\sc Proof.} We may, even if it means translating $u_\e$, assume the normalisation condition
\begin{equation}
\label{e3.400}
\min_{[-L,0]}v_\e(0,y)=\e.
\end{equation}
Therefore, $f_\e(u)\equiv0$ on $\RR_-\times[-L,0]$, and $u_\e$ solves a purely linear equation. If
$$
\rho\leq\min(\frac1D,\frac1d)c,
$$
then 
$$(\underline u(x),\underline v(x,y)):=\biggl(1-(1-\e)e^{\rho x}\biggl)(\frac1\mu,1)
$$ is a subsolution to \eqref{e1.2} on $\RR_-$, thus $u(x,y)\geq \underline u(x)$ on $\RR_-$. Choosing $\rho$ the be the above minimum, we obtain
$$
\partial_xv(0,y_e)\leq-(1-\e)\min(\frac1D,\frac1d)c,
$$
where $y_\e$ is a point in $[-L,0]$ where the minimum in \eqref{e3.400} is attained.

\noindent On the right half of $\Sigma$, another simple sub-solution to \eqref{e1.2} is 
$$
(\underline u_\e(x),\underline v_\e(x,y))=\e p(\frac{x}\e)(\frac1\mu,1),
$$
with
$$
-p''+p\varphi(p)=0\ \hbox{for $\xi>0$,}\quad p(0)=1,\ p(+\infty)=0.
$$
Because it is convex, it is a sub-solution to the equation for $v$, and because it is decreasing it is also a sub-solution of the equation inside 
the right half of $\Sigma$. The exchange condition for $v$ is automatically satisfied. The classical sliding argument (slide $(u_\e,v_\e)$ until  the two components exceeds $(\underline u_\e,v_\e)$, then slide back until one of the two components reaches a contact point) yields $(u_\e,v_\e)\geq(\underline u_\e,v_\e)$. Thus
$$
\partial_x v_\e(0,y_\e)\geq \underline v_\e'(0)=p'(0)=-\sqrt 2.
$$
This implies
$$
(1-\e)\min(\frac1D,\frac1d)c\leq\sqrt2,
$$
the sought for bound. \hfill$\square$

\noindent The next step is a uniform gradient bound on $v_\e$. Notice that $u_\e$, that satisfies a linear ODE with bounded right handside, thus is uniformly $C^{1,1}$ at this stage. 
\begin{prop}
\label{p3.2000}
There is $M>0$ universal such that $\vert\nabla v_\e\vert\leq M$ in $\overline\Sigma$. 
\end{prop}
In the sequel we will, for convenience, drop the subscript ${}_\e$ for $c_\e$, $u_\e$ and $v_\e$.
The first ingredient is a gradient bound in $\{v\leq\e\}$, away from the road.

\begin{lm}
\label{l3.1000} 
Consider $\lambda\in(0,1)$ and a point $(x_0,x_0)\in\overline\Sigma$ such that 
$
u(x_0,y_0)=\lambda\e.
$
Assume that $y_0\leq-2\e.$  Then we have
$$
0\leq u(x,y)\leq C\e,\quad (x,y)\in B_\e(x_0,y_0),
$$
and
$$
\vert\nabla u(x_0,y_0)\vert\leq C,
$$
for some universal $C>0$.
\end{lm}
\noindent{\sc Proof. } Notice, following \cite{BCN}, that the case $y_0$ close to $-L$ is not really an issue, because one may, thanks to the Neumann condition at $y=-L$, extend the function $u$ evenly in $y$ to $\RR\times (-2L,0)$. This said, we do the classical Lipschitz scaling
\begin{equation}
\label{e3.1700}
v(x_0+\e\xi,\e\zeta)=\e p(\xi,\zeta),\quad (\xi,\zeta)\in B_2(0),
\end{equation}
 and $p$ solves
$$
-d\Delta p+\e c\partial_\xi p+p\varphi(p)=0\quad(\xi,\zeta)\in B_2(0),\ p(0,0)=\lambda.
$$
Then (\cite{BCN} again), from the Harnack inequality, $p$ is universally controlled in $B_1(0)$, hence $\nabla u(x_0,y_0)=\nabla p(0,0)$ is universally controlled. \hfill$\square$
\begin{lm}
\label{l3.200} 
Consider $\lambda\in(0,1)$, $y_0\in[-\e,0]$ and $x_0\in\RR$ such that 
$$
u(x_0,y_0)=\lambda\e.
$$
Then we have
$$
0\leq u(x,y)\leq C\e\lambda,\quad (x,y)\in B_\e(x_0,y_0),\ y\leq0,
$$
for some $C>0$ that depends neither on $\e$, nor $\lambda$.
\end{lm}
\noindent{\sc Proof.} Recall that $u$ and $v$ are bounded independently of $\e$. We redo the scaling \eqref{e3.1700}, leaving $u$ untouched. The only thing that has to be examined is the Neumann condition for $p$, which reads
\begin{equation}
\label{e3.1701}
p_\zeta(\xi,0)+\e p(\xi,\zeta)=u(x_0+\e\xi,0).
\end{equation}
We make the slight abuse of notations consisting in denoting by $u(\xi,\zeta)$ the function $u(x_0+\e\xi,y_0+\e\zeta)$, this function is clearly $C^{1,\alpha}$ uniformly in $\e$. We may indeed subtract from $p$ any suitable harmonic function $V$ satisfying the Neumann condition \eqref{e3.1701}, thus $V$ is uniformly $C^{2,\alpha}$
in $B_{3/4}$. We then apply the argument of Lemma \ref{l3.200} to $v-V$. \hfill$\square$

\noindent These two lemmas lead to an effortless

\noindent{\bf Proof of Proposition \ref{p3.2000}.} Let $\Gamma_\e$ be the curve $\{v=\e\}$. The function $v_y$ is bounded on $\Gamma_\e$, as well as on $\partial\Sigma\backslash\Gamma_\e$, moreover it satisfies $-\Delta y_y+c\partial_xv_y=0$ to the left of $\Gamma_\e$, thus it is bounded everywhere. The same argument applies to $v_x$, but for this quantity we use the Robin condition $d\partial_y v_x+v_x=u_x$ (recall that $u_x$ is uniformly bounded) and the Neumann conditon for $v_x$ at the bottom of $\Sigma$. \hfill$\square$

\subsection{Convergence of the approximating sequence} 
The last ingredient that we need is a uniform lower bound on the travelling vave velocity. In this proposition we put the subscripts ${}_\e$ back in.
\begin{prop}
\label{p3.3000}
There is $c_0>0$ such that $c_\e\geq c_0$.
\end{prop} 
\noindent{\sc Proof.} We start with the following identity
\begin{equation}
\label{e3.40}
c_\e=\frac1{L+d\mu}\int_\Sigma f_\e(v_\e)dxdy,
\end{equation}
obtained by integrating the system for $(u_\e,v_\e)$ over $\Sigma$.
We may normalise $v_\e$ so that
$$
v_\e(0,-\frac{L}4)=\e.
$$
From \cite{BCN}, the family of measures
$$
\sigma_\e=\un_{B_{L/4}}(0,-\frac{L}4)f_\e(v_\e)dxdy
$$
converges (possibly up to a subsequence), in the measure sense, to the image of the Lebesgue measure on a locally BV graph
$\{h(y),y)\}$. Thus the sequence $\biggl((L+d\mu)c_\e\biggl)_\e$ converges, still up to a subsequence, to a limit that is larger than the length of $\Gamma$ inside $B_{L/4}(0,-\di\frac{L}4)$, a positive number. \hfill$\square$

Putting everything together, we may give the

\noindent{\sc Proof of Theorem \ref{t2.1}.} Possibly up to a subsequence, the sequence $(u_\e,v_\e)_\e$ converges, uniformly in $\overline\Sigma$, as well as in $H^1_{loc}(\Sigma)$ weakly, to a function $(u(x),v(x,y))$ that is both Lipschitz and in $H^1_{loc}(\Sigma)$. Notice that $u$ is much smoother, it is $C^{2,1}$. Let us repeat here that the family of measures $(f_\e(v_\e)dxdy)_\e$ converges (\cite{BCN} again), in every set of the form $\RR\times K$ ($K$ a compact subset of $[-L,0)$), to the length measure of a graph
$$
\Gamma=\{(h(y),y)\,\ -L\leq y<0\}.
$$
For every $\delta\in(0,\di\frac{L}4)$, identity \eqref{e3.40} implies
\begin{equation}
\label{e3.41}
\int_{\RR\times[-L,-\delta]}f_\e(u)dxdy\leq (L+d\mu)c_\e.
\end{equation}
Passing to the limit $\e\to0$ yields
\begin{equation}
\label{e3.42}
\int_{-L}^{-\delta} \sqrt{1+(h'(y))^2}dy\leq (L+d\mu)c.
\end{equation}
Thus the function $h$ belongs to $B\!V([-L,0))$, thus can be extended by continuity to $y=0$. We may always assume
$
h(0)=0.
$ 
It remains to prove that $\partial\{u>0\}$ is analytic. Let us set 
$$
F_\e(v)=\int_0^vf_\e(v)dv,
$$
we claim that, for every $\e>0$, $v_\e$ is a local minimiser of the energy
$$
\int e^{-c_\e x}\biggl(\frac12\vert \nabla v\vert^2+F_\e(v)\biggl)dxdy.
$$
More precisely, for every ball $B$ whose closure is included in $\Sigma$, then $v_\e$ minimises the energy
\begin{equation}
\label{e2.2000}
J_\e(\phi,B)=\int_Be^{-c_\e x}\biggl(\frac12\vert\nabla \phi\vert^2+F_\e(\phi)\biggl),
\end{equation}
over all functions $\phi\in H^1(B)$ whose trace on $\partial B$ is $v_\e$. 
This is easily seen from the monotonicity of $v_\e$ in $x$, and a sliding argument. On the other hand we use the two following facts, taken from \cite{BCN}: first,
the family $(v_\e)_\e$ is compact in $H^1_{loc}(\Sigma)$, (ii). there is a uniform nondegeneracy property:
$$
v_\e(x,y)\geq kd\biggl((x,y),\{v_\e\in(a\e,b\e)\}\biggl).
$$
This implies that, for every ball $B$ inside $\Sigma$ we have:
$$
\lim_{\e\to 0}J_\e(v_\e,B)=\int_B\frac{e^{-cx}}2\vert\nabla v\vert^2dxdy+\int_{\{v>0\}\cap B}e^{-cx}dxdy:=J(v,B).
$$
On the other hand, for every $\phi\in H^1(B)$ whose trace on $B$ is $v$, consider $\phi_\e\in H^1(B)$, whose trace on $B$ is $v_\e$, and
such that the family $(\phi_\e)_\e$ converges to $\phi$ in $H^1(B)$. We have $J_\e(v_\e,B)\leq J_\e(\phi,B)$, which implies, sending $\e\to0$: 
$J(v,B)\leq J(\phi,B)$. Thus, $v$ is a local minimiser of $J$. However, $J$ is of the type of functionals treated in \cite{ACF}: its results are applicable, 
which imples the analyticity of $\partial\{v>0\}$ inside $\Sigma$.

\noindent Finally, because $v$ is Lipschitz, and because of the uniform convergence of $(v_\e)_\e$ to $v$, we have $v(x,0)=0$ if $x>0$. 
\hfill$\square$

\noindent{\bf Remark.} {\it The above argument also explains the loss of the exchange condition for $v$ to the right of the free boundary, simply because the free boundary is forced to hit the road. There is a boundary layer in which the condition $v=0$, located on a curve very close to the road, eventually overcomes the exchange condition in the limit $\e\to0$.}
\section{The two species model: the free boundary near (0,0)}
The main feature of Model \eqref{e1.200} is that the equation inside $\Sigma$, together with the free boundary conditions, can be studied in the neighbourhood of a free boundary point up to the top of $\Sigma$ {\it via} Lipschitz rescalings: $(x,y)=\delta(\xi,\zeta)$. This will enable us to show, in a relatively easy way, the linear behaviour of the free boundary in the vicinity of the road. Let us first state a rigidity result in 2 space dimensions.
\begin{thm}
\label{t4.52}
Pick $\lambda\in[0,1)$. Let $u(x,y)$ solve
\begin{equation}
\label{e4.501}
\begin{array}{rll}
\Delta u=&0\quad (x,y)\in \RR\times\RR_-\cap\{u>0\}\\
\vert\nabla u\vert=&1\quad (x,y)\in\partial\{u>0\},\\
u_y(x,0)=&\lambda\quad \hbox{if $x<0$},\\
u(x,0)=&0\quad\hbox{if $x>0$.}
\end{array}
\end{equation} 
Also assume that $\partial_xu\leq0$. Then we have:
$$
\partial\{u>0\}=\{y=-\frac{\sqrt{1-\lambda^2}}{\lambda}x,\ x<0\},
$$
and 
$$
u(x,y)=\biggl(\sqrt{1-\lambda^2}x+\lambda y\biggl)^+.
$$
\end{thm}
\noindent{\bf Proof.} Let $\Omega$ be the positivity set of $u$, identify $\RR^2$ with the complex plane and set
$$
f(z)=u_y(x,y)+iu_x(x,y),\quad z=x+iy.
$$
Then $f$ is analytic in $\Omega$, let us assume that it is a nonconstant function. It is open from $\Omega$ onto its image, and  maps $\Gamma$ onto (a portion of) the unit circle, that we call $\gamma_1$, whereas it maps the negative $x$ axis onto (a portion of) the vertical line $\{{\mathrm{Im}}Z=\lambda\}$, that we call $\gamma_2$. From its connectedness,  $\gamma_1\cup\gamma_2$ (but not only) enclose $f(\Omega)$. And, because $\partial_xu\leq 0$,  $f(\Omega)$ is bounded by $\gamma_1\cup\gamma_2\cup\gamma_3$, where $\gamma_3$ is a non void, possibly very irregular curve. One of its end points is on $\gamma_1$ - call it $z_1$; thus $\gamma_1$ is the segment $[z_1,\sqrt{1-\lambda^2}+i\lambda]$ - and the other on $\gamma_2$ - call it $z_2$. It stays in the upper half of the complex plane, and is constructed as follows: for each  ray $D$ starting from $\sqrt{1-\lambda^2}+i\lambda$, and not meeting $\gamma_1\cup\gamma_2$, let $Z_D$ be the furthest point from $\sqrt{1-\lambda^2}+i\lambda$ in $\overline{f(\Omega)}$: the set of all such $Z_D$'s determines $\gamma_3$. Two cases are to be considered.

\noindent{\bf Case 1.} {\it $\gamma_3$ has points strictly above the horizontal axis.} Let $D$ be a ray starting from $\sqrt{1-\lambda^2}+i\lambda$, and let $Z_0$ be the furthest point in $D\cap\gamma_3$. If $Z_0$ has a pre-image by $f$, we have a contradiction because $f$ is open. If $Z_0$ has no pre-image, this means the existence of a sequence $(z_n)_n$ such that
$$
\lim_{n\to+\infty}\vert z_n\vert=0\ \hbox{or}\ \lim_{n\to+\infty}\vert z_n\vert=+\infty,
$$
and such that 
$\di\lim_{n\to+\infty}f(z_n)=Z_0.
$
The sequence of locally uniformly bounded functions
\begin{equation}
\label{e3.1000}
f_n(z)=f(\frac{z}{\vert z_n\vert})\quad\hbox {if $\vert z_n\vert\to0$},
\end{equation}
or
\begin{equation}
\label{e3.1001}
f_n(z)=f({z}{\vert z_n\vert})\quad\hbox{if $\vert z_n\vert\to+\infty$},
\end{equation}
converges, up to a subsequence, uniformly locally to an analytic function $f_\infty(z)$, by Montel's theorem. The function $f_\infty$ takes the value $Z_0$, but it is also nonconstant analytic, in a possibly different positivity set $\Omega_\infty$. However, by definition, $Z_0$ is still the furthest point between $D$ and $f_\infty(\Omega_\infty)$, and the openness of $f_\infty$ yields the contradiction. Thus  $f$ is constant, and so is $\nabla u$.

\noindent{\bf Case 2.} {\it $\gamma_3$ meets the horizontal axis}. Let $\mu\in[0,1]$ belong to $\gamma_3$. Similarly as above, there is a sequence $(z_n)_n$ whose modulus goes to 0 or $+\infty$, such that
$
\di\lim_{n\to+\infty}f(z_n)=\mu.
$
Let $f_n$ be defined as in \eqref{e3.1000} or \eqref{e3.1001}, and let $f_\infty$ be one of its limits. From nondegeneracy (i.e. linear growth from each point of the free boundary, which is, in this case, inherited from its analyticity) there is a nontrivial solution $u_\infty$ of \eqref{e4.501}, and a point $(x_0,y_0)$ in the positivity set of $u_\infty$ such that 
$$
\partial_xu_\infty(x_0,y_0)=0.
$$
Because $\partial_xu_\infty\leq0$, it is zero everywhere. The only possibility for $\Gamma_\infty$ is that it is horizontal, a contradiction to the fact that it has to meet the $x$ axis.
\hfill $\square$

In particular, we have, for $\lambda=0$:
$$
\partial\{u>0\}=\{y=0\},\quad
u(x,y)=y^-.
$$
The same argument as in Theorem \ref{t4.52} yields the following
\begin{cor}
\label{c3.10}
Pick $\lambda\in[0,1]$, and consider a  solution $(\varphi(y),u(x,y))$, with $\varphi(y)<0$ and $\partial_xu\leq0$, of
\begin{equation}
\label{e4.601}
\begin{array}{rll}
\Delta u=&0,\ u>0\quad (x<\varphi(y))\\
\vert\nabla u\vert=&1\quad (x=\varphi(y)),\\
u_y(x,0)=&\lambda\quad (x\in\RR).\\
\end{array}
\end{equation} 
Then $\varphi$ is constant negative, and $\lambda=1$.
\end{cor} 
The consequence of these rigidity properties for 2D global solutions is the

\noindent{\bf Proof of Theorem \ref{t2.2}.} Let $(c,u,v,\Gamma)$ the solution of \eqref{e1.200}. We claim that $\Gamma$ hits the road at a point that we may assume, by translational invariance, to be $(0,0)$. Let us prove that 
$u(x_0)>0$. Indeed, assume the contrary and that $u(x_0)=0$. Consider a Lipschitz blow up of $v$ around $(x_0,0)$:
\begin{equation}
\label{e4.2040}
v_\delta(\xi,\zeta)=\frac{v(\delta\xi,\delta\zeta)}{\delta}.
\end{equation}
Call $\Gamma_\delta$ its free boundary, let us first show that, as $\delta\to0$, it does not collapse on the $x$ axis. In other words, we want to show the existence of $\rho_0>0$ such that, if $\Gamma_\delta$ is the positivity set of $u_\delta$ within $B_1(0)$, then 
$$
B_{\rho_0}(-1,-2\rho_0)\subset\Omega_\delta.
$$
Assume that it is not the case, then we have
$
\di\lim_{\delta\to0}v_\delta(\xi,0)=0$,
uniformly in every compact  in $\xi$, just because $v$ is Lipschitz. And so, by the exchange condition, we have
$
\di\lim_{\delta\to0}\partial_yv_\delta(-1,0)=0.
$
On the other hand, let $\zeta_\delta$ be the smallest $\zeta<0$ such that $(-1,\zeta)\in\Gamma_\delta$. Clearly, $\zeta_\delta$ tends to 0 as $\delta\to0$. Rescale $v_\delta$ with $\vert\zeta_\delta\vert$, taking $(-1,\zeta_\delta)$ as the origin. That is, we set
$$
v_\delta(\xi,\zeta)=\vert\zeta_\delta\vert V_\delta\biggl(\frac{\xi+1}{\vert\zeta_\delta\vert},\frac{\zeta-\zeta_\delta}{\vert\zeta_\delta\vert}\biggl).
$$
(A subsequence of) the sequence $(V_\delta)_\delta$ will converge, as $\delta\to0$, to a solution $V_\infty$ of \eqref{e4.601}, with $\lambda=0$, something that Corollary \ref{c3.10} excludes.
As $\delta\to0$, (a subsequence of) $(v_\delta)_\delta$ converges locally uniformly to a global solution $V_\infty$ of \eqref{e4.501}, with $\lambda=0$. By Theorem \ref{t4.52}, we have 
$
V_\infty(\xi,\zeta)=\xi^-.
$
Thus, the equation for $u$ in a neighbourhood of $x=x_0$ is
$$
\begin{array}{rll}
-Du''+cu'+\mu u\sim& x_0-x\quad x<0\\
u(0)=u'(0)=&0.
\end{array}
$$
This implies
$
u(x)\sim \di\frac{x^3}{6D}$ for $x<0$, close to $0$.
Thus $u(x)<0$ in a neighbourhood of $0$, a contradiction with the positivity of $u$. So, $u(0)>0$, let us set
$$
\lambda=\mu u(0)>0.
$$
Let us come back to the blow up \eqref{e4.2040}. Each converging subsequence as a solution of \eqref{e4.501} as a limit, which is unique 
by virtue of Theorem \ref{t4.52}. So, the whole blow-up converges to 
$$
V_\infty(x)=\biggl(\lambda y-\sqrt{1-\lambda^2}x\biggl)^+,
$$
which implies, coming back to the solution $(c,u,v,\Gamma)$ of \eqref{e1.200}, that $\Gamma$ has slope 
$
\gamma=\di\frac{\sqrt{1-\lambda^2}}\lambda
$
at $(0,0)$. This proves Theorem \ref{t2.2}. \hfill$\square$
\section{The travelling wave of the single species model}
We will, as in Section 2, construct a solution to \eqref{e1.2} as the limit, as $\e\to0$, of a sequence $(c_e,u_\e)_\e$ of solutions of 
\begin{equation}
\label{e3.2}
\left\{
\begin{array}{rll}
-d\Delta u+c\partial_x u+f_\e(u)=&0\quad(x,y)\in\Sigma\\
\ \\
-Du_{xx}+c\partial_xu+1/\mu u_y=&0\quad\hbox{for $x\in\RR$, $y=0$}\\
u_y(x,-L)=&0\\
\ \\
u(-\infty,y)=&1,\quad u(+\infty,y)=0.
\end{array}
\right.
\end{equation}
the function $f_\e$ being as in \eqref{e3.1}. For every $\e>0$, \eqref{e3.2} has (Dietrich \cite{Di1} once again) a unique solution $(c_\e,u_\e)$ such that $c_\e>0$, $0<u_\e<1$ and $\partial_xu_\e<0$; we will show that, up to a subsequence, $(c_\e,u_\e)$ converges, as $\e\to0$, to a solution of \eqref{e1.2}. Our first task is to show that $c_\e$ is uniformly bounded from above, that is:
\begin{prop}
\label{p3.1}
There is $K>0$ independent of $\e$ such that $c_\e\leq K$.
\end{prop}
\noindent{\sc Proof.} Similar to that of Proposition \ref{p3.1000}. In $\RR_-\times[-L,0]$, we work with the sub-solution
$$
\underline u(x):=1-(1-\e)e^{\rho x}
$$
with
$$
\rho\leq\min(\frac1D,\frac1d)c,
$$
 whereas, to the right of $\Sigma$, we work with
$$
\underline u_\e(x)=\e p(\frac{x}\e),
$$
with
$$
-p''+p\varphi(p)=0\ \hbox{for $\xi>0$,}\quad p(0)=1,\ p(+\infty)=0.
$$
Because it is convex, it is a sub-solution to the Wentzell boundary condition on $\RR_+$, and because it is decreasing it is also a sub-solution of the equation inside 
the right half of $\Sigma$. \hfill$\square$

\noindent The next step is a uniform gradient bound on $u_\e$:
\begin{thm}
\label{t4.1}
Assume $D\geq d$. There is $M>0$ universal such that $\Vert\nabla u_\e\Vert_\infty\leq M$.
\end{thm}
This does not result from the application of known theorems, and will be rather involved. So, we first assume that it holds, and finish the construction of the travelling wave of \eqref{e1.2}.
\subsection{Construction of a solution to the one species model, given Theorem \ref{t4.1}}
\noindent Let us first state the equivalent of Proposition \ref{p3.3000}
\begin{prop}
\label{p3.3}
There is $c_0>0$ such that $c_\e\geq c_0$.
\end{prop} 

\noindent{\sc Proof of Theorem \ref{t2.3}.} Arguing as in the proof of Theorem \ref{t2.1}, we obtain a solution $u$ to \eqref{e1.2}, at least for the free boundary problem inside, with a free boundary $(h(y),y)$ that meets the top of $\Sigma$ at a point that we may assume to be 0. So, we have $u(x,0)=0$ if $x>0$.

\noindent To the left of $\Gamma\cap\Sigma$ we have $u>0$. We infer that $u(x,0)>0$ if $x<0$. Indeed, assume $u(x,0)=0$ for $x<\bar x<0$, and $u(x,0)>0$ if $x>\bar x$. As a consequence of the definition of $\Gamma$, there is $r\in(0,\vert\bar x\vert)$ such that 
$$
u>0\ \hbox{in}\ B_r(\bar x-r,0)\cap\{y<0\}.
$$
Let us come back to $u_\e$, the situation is the following:
\begin{equation}
\begin{array}{rll}
-d\Delta u_\e+c_\e\partial_xu_\e=&0\quad \bar x-r\leq x\leq\bar x,\ -r<y<0\\
-D\partial_{xx}u_\e+c_\e\partial_xu_\e+\partial_yu_\e=&0,\quad \bar x-r\leq x\leq\bar x\\
\ \\
u_\e(x,0)=&o_{\e\to0}(1),\quad \bar x-r\leq x\leq\bar x\\
u_\e(x,y)\geq \lambda_r\vert y\vert\ \hbox{for some $\lambda_r>0$,}
\end{array}
\end{equation}
the last inequality being valid because of the Hopf Lemma for $u$ applied at every point of the road between $(\bar x-r,0)$ and $(\bar x,0)$ on the one hand, and the fact, on the other hand, that the convergence of $u_\e$ to $u$ is better than uniform in $\{\bar x-r\leq x\leq\bar x,\ -r\leq y\leq0\}$ - the Schauder estimates being valid in this area. As a consequence we have
\begin{equation}
\label{e3.43}
\partial_y u_\e(x,0)\leq-\lambda_r\ \hbox{for}\ \bar x-r\leq x\leq\bar x,\ -r\leq y<0.
\end{equation}
Because $u_x(\bar x-r,0)=0$ we have 
$u_\e(\bar x-r,0)=o_{\e\to0}(1).
$
Integration of the Wentzell condition for $u_\e$ between $\bar x-r$ and $\bar x$, and taking \eqref{e3.43} into account we obtain, when $\e>0$ is small enough:
$$
u_\e(\bar x,0)<0,
$$
a contradiction. 

\noindent Let us finally show that $u$ solves the full free boundary problem \eqref{e1.2}. From identity \eqref{e3.40} and the boundedness of the sequence $(c_\e)_\e$ we obtain that $(f_\e(u_\e)dxdy)_\e$ converges (still up to a subsequence) to the length measure on $\Gamma$, plus a possible finite measure on the boundary $\RR\times\{y=0\}$. We note that $u(x,0)$ is at least $C^{1,1}$: indeed, $\partial_yu_\e$ is uniformly bounded, due to Theorem \ref{t4.1}. Thus $\partial_{xx}u_\e(.,0)$ is also uniformly bounded, a property that passes on to $\partial_{xx}u$. So, the contribution of the limiting measure on $\RR\times\{y=0\}$ is nonexistent, which proves that the Wentzell condition is satisfied both a.e. and in the distributional sense. The rest of the system is proved to be solved just as in Theorem \ref{t2.1}. \hfill$\square$
\subsection{Gradient bound in the region $u\sim\e$}
The scheme of the proof is still, roughly, that of the Lipschitz bound in \cite{BCN}. However, here, the vicinity of the road requires a special treatment and this is where we will, eventually, need $D\geq d$.  
In the sequel we will, for convenience, drop the subscript ${}_\e$ for $c_\e$ and $u_\e$.
First, recall the gradient bound in $\{u\leq\e\}$, away from the road.

\begin{lm}
\label{l3.1} 
Consider $\lambda\in(0,1)$ and a point $(x_0,y_0)\in\overline\Sigma$ such that 
$$
u(x_0,y_0)=\lambda\e.
$$
Assume that $y_0\leq-2\e.$  Then we have
$$
0\leq u(x,y)\leq C\e,\quad (x,y)\in B_\e(x_0,y_0),
$$
and
$
\vert\nabla u(x_0,y_0)\vert\leq C,
$
for some universal $C>0$.
\end{lm}
\noindent{\sc Proof. } Set 
\begin{equation}
\label{e3.17}
p(\xi,\zeta)=\e u(x_0+\e\xi,y_0+\e\zeta),
\end{equation}
and apply \cite{BCN}.
 \hfill$\square$

\noindent The main part of the task is therefore to bound $\nabla u$ at distance less than $\e$ from the road. Let us start with the most extreme case, i.e. a point on the road. 
\begin{lm}
\label{l3.2} 
Consider $\lambda\in(0,1)$ and $x_0\in\RR$ such that 
$$
u(x_0,0)=\lambda\e.
$$
Then we have
$$
0\leq u(x,y)\leq C\e\lambda,\quad (x,y)\in B_\e(x_0,y_0),\ y\leq0,
$$
for some $C>0$ that depends neither on $\e$, nor $\lambda$.
\end{lm}
\noindent{\sc Proof.} This looks, at first sight,  like an innocent repetition of the previous lemma. However one quickly realises that the Lipschitz scaling yields the boundary condition
$$
-Dp_{\xi\xi}+\e cp_\xi+\e/\mu p_\zeta=0,
$$
so that half of the Wentzell condition is lost when one sets $\e=0$. Of course this yields that $p$ is linear on the road, hence constant in order to keep its positivity. However it does not prevent large gradients in $y$: assume for definiteness that $\lambda=1$;  if $p(-1,0)=M>>1$, then 
$$
q(\xi,\zeta)=\frac{p(\xi,\zeta)}M,
$$
and $q$ solves 
$$
-\Delta q+q\varphi(Mq)=0
$$
inside the cylinder, while keeping the condition $q_{\xi\xi}=0$. So, it is an asymptotic global solution for infinite $M$, which implies that one cannot hope for a uniform bound for $p$ in this setting.

\noindent So, this time we use the mixed Lipschitz-quadratic scaling
\begin{equation}
\label{e3.117}
u(x_0+\sqrt\e\xi,\e\zeta)=\e \lambda p(\xi,\zeta),\quad (\xi,\zeta)\in B_2(0),\ \zeta<0
\end{equation}
 and $p$ solves
\begin{equation}
\label{e3.8}
\left\{
\begin{array}{rll}
-d\e p_{\xi\xi}-dp_{\zeta\zeta}+\sqrt\e c\partial_\xi p+p\varphi(\frac{p}\lambda)=&0\quad(\xi,\zeta)\in \cap\{\zeta<0\}\\
-Dp_{\xi\xi}+\sqrt\e c\partial_\xi p+1/\mu p_\zeta=&0\quad\hbox{for $-2<\xi<2$, $\zeta=0$}\\
p(0,0)=&\lambda.
\end{array}
\right.
\end{equation}
We claim that $p(0,-1)$ is universally bounded, both with respect to $\lambda$ and $\e$. Call $M_{\lambda,\e}$ this quantity, we have 
$$
p(\xi,-1)\geq M_{\lambda,\e}\ \hbox{for $\xi\leq0$.}
$$
Assume that a subsequence of $(M_{\lambda,\e})_{\lambda,\e}$ (that we still relabel $(M_{\lambda,\e})_{\lambda,\e}$) with $\e\to0$, grows to infinity. We claim the existence of $\xi_{\lambda,\e}<0$, going to 0 as $\e\to0$, such that
\begin{equation}
\label{e3.10}
p(\xi_{\lambda,\e},)\geq kM_{\lambda,\e},
\end{equation}
for some universal $k>0$. If such is not the case, then there is $\xi_0<0$ universal such that
$$
\lim_{\e\to0}\frac{p(\xi,0)}{M_{\lambda,\e}}=0,\  \hbox{uniformly in $\xi\in[\xi_0,0]$.}
$$
This entails the existence of  a constant $\gamma_0>0$, universal, such that
\begin{equation}
\label{e3.9}
p_\zeta(\xi,0)\leq-\gamma_0 M_{\lambda,\e},\ \hbox{$\xi\in[\xi_0+\sqrt\e,-\sqrt\e].$}
\end{equation}
Indeed, we have $p(\xi,\zeta)\geq\underline p(\xi,\zeta)$, where 
$$
\begin{array}{rll}
-d\e\underline p_{\xi\xi}-d\underline p_{\zeta\zeta}+\sqrt\e^{3/2} c\underline p_\xi+\underline p=&0\quad (\xi\in(\xi_0,0)\times(-1,0))\\
\underline p(\xi,0)=p(\xi,0),\ \underline p(\xi,-1)=&M_{\lambda_\e}\quad (\xi\in(\xi_0,0))\\
\underline p(0,\zeta)=p(\xi_0,\zeta)=&0\quad (\zeta\in(-1,0)).
\end{array}
$$
Rescaling in $\xi$ - so as to recover a fully elliptic equation for $p$ - yields that $\underline p_{\zeta}(\xi,0)$ satisfies an estimate of the type \eqref{e3.9}.

\noindent Now, from Rolle's theorem, there is $\xi_{\lambda,\e}'\in (\di\frac{\xi_0}2,\di\frac{\xi_0}4)$ such that $p_\xi(\xi_{\lambda,\e}',0)=o(M_{\lambda,\e})$ (recall that $p$ is an $o(M_{\lambda,\e})$ on $(\xi_0+\e,-\e)$). So we have, on $(\di\frac{3\xi_0}4,\xi_{\lambda,\e}')$:
$$
Dp_{\xi\xi}(\xi,0)-\e cp_\xi(\xi,0)\leq-\frac{\gamma_0M_{\lambda,\e}}{\mu},\quad p_\xi(\xi_{\e,\lambda}',0)=o(M_{\lambda,\e}),\ p(\xi_{\e,\lambda}'(\xi_{\e,\lambda}',0)=o(M_{\lambda,\e}).
$$
Integration of this very simple differential inequality between $\xi_{\e,\lambda}'$ and $\di\frac{3\xi_0}4$ yields
$$
p(\di\frac{3\xi_0}4,0)=o(M_{\lambda,\e})-\frac12\gamma_0M_{\e,\lambda}(\frac{3\xi_0}4-\xi_{\lambda,\e}')^2\leq O(1)-\frac{9\xi_0^2}{16}<0,
$$
a contradiction.

\noindent So, we have found $\xi_{\e,\lambda}$, going to 0 as $\e\to0$, such that \eqref{e3.10} holds. Thus there exists $\xi_{\e,\lambda}''\in(\xi_{\e,\lambda}',0)$ such that
$$
p(\xi_{\e,\lambda}'',0)\in(1,2),\quad\lim_{\e\to0}p_\xi(\xi_{\e,\lambda}''0):=a_{\lambda,\e}=-\infty.
$$
Now we notice that $p_\zeta(\xi,0)$ is less than some universal $C>0$, for $\xi\in(\xi_{\e,\lambda}'',+\infty)$: indeed, set 
$$
\tilde p(\xi)=\left\{
\begin{array}{rll}
p(\xi_{\e,\lambda}'',0)\ &\hbox{if $\xi\leq\xi_{\e,\lambda}''$,}\\
p(\xi,0)\ &\hbox{if $\xi\geq\xi_{\e,\lambda}''$,}
\end{array}
\right.
$$
so that $\tilde p(\xi)\leq C$ universal on $(\xi_{\e,\lambda}'',1)$. We have 
$
p(\xi,\zeta)\leq \underline q(\xi,\zeta)
$ where 
$$
\begin{array}{rll}
-d\e\underline q_{\xi\xi}-d\underline q_{\zeta\zeta}+\sqrt\e^{3/2} c\underline q_\xi+\underline q=&0\quad (\xi\in(-1,1)\times(-1,0))\\
\underline q(\xi,0)=\tilde p(\xi),\ \underline p(\xi,-1)=&0\quad (\xi\in(-1,1))\\
\underline p(\pm 1,\zeta)=&0\quad (\zeta\in(-1,0)),
\end{array}
$$
which satisfies 
$$\underline q_\zeta(\xi,0)\leq C,\ \hbox{$C>0$ universal.}
$$ 
So we have this time
$$
Dp_{\xi\xi}(\xi,0)-\e cp_\xi(\xi,0)\leq-C\ \hbox{on $(\xi_{\e,\lambda}'',1)$},
$$
while
$$
p_\xi(\xi_{\e,\lambda}'',0)\in (0,2],\ p(\xi_{\e,\lambda}'(\xi_{\e,\lambda}',0)=a_{\lambda,\e}\to-\infty.
$$
Again, integration of the differential inequality on $(\xi_{\e,\lambda}'',1)$ yields the existence of $\xi_{\e,\lambda}'''\geq \xi_{\e,\lambda}''$
such that $p(\xi_{\e,\lambda}''',)<0$, a contradiction. Thus our claim that $p(0,-1)$ is universally bounded is proved.

\noindent It remains to see that $p(0,\zeta)$ is universally bounded for $\zeta\in[-1,0]$. First, notice that we have in fact proved that
$$
p(\xi,-1)\ \hbox{is universally bounded for $-1\leq\xi\leq1$.}
$$
We have also proved that 
$$
p_\zeta(\xi,-1)\leq C,\ \hbox{$C>0$ universal.}
$$
This remark in hand we may repeat the above argument, replacing $-1$ by $\zeta$. \hfill$\square$

\noindent Instructed by Lemma \ref{l3.2}, we may now bound $u$ in an $\e$-vicinity of the road.
\begin{lm}
\label{l3.3} 
Consider $\lambda\in(0,1)$ and $(x_0,y_0)\in\RR\times(-\e,0)$ such that 
$$
u(x_0,y_0)=\lambda\e.
$$
Then we have
$$
0\leq u(x,y)\leq C\e\lambda,\quad (x,y)\in B_{2\e}(x_0,y_0),\ y\leq0,
$$
for some $C>0$ that depends neither on $\e$, nor on $\lambda$.
\end{lm}
\noindent{\sc Proof.} We start again the mixed scaling \eqref{e3.117}, so that
$p$ solves \eqref{e3.8}. Set
$$
y_0=\e\zeta_0,\ \hbox{thus}\ u(x_0,y_0)=\e p(0,\zeta_0)=\e\lambda.
$$
The only thing that we have to prove is that $p(0,0)$ is universally bounded from above and below. The case
$$
\lim_{\e\to0}\frac{p(0,0)}{p(0,\zeta_0)}=0
$$ 
amounts to the previous lemma, so it only remains to exclude the case
\begin{equation}
\label{e3.21}
\lim_{\e\to0}\frac{p(0,0)}{p(0,\zeta_0)}=+\infty.
\end{equation}{
This time it will be more useful to work in the Lipschitz scaling \eqref{e3.17}. Set
$$
M_{\e,\lambda}=\frac{p(0,0)}{p(0,\zeta_0)},\quad q(\xi,\zeta)=\frac{p(\xi,\zeta)}{M_{\e,\lambda}}, 
$$
so that $q$ solves
$$
-d\Delta q+\e c\partial_\xi q+q\varphi(M_{\e,\lambda}q)=0\quad(\xi,\zeta)\in B_2(0),\ q(0,0)=1.
$$
Now, notice the very simple bound for $\nabla q$: The supremum of $f_\e$ being of order $\e$, the elliptic estimates \cite{LTru} for the 
original Wentzell problem  
\eqref{e3.2} yield
$$
\Vert D^2 u\Vert_{W^{2,p}(B\cap\bar\Sigma)}\leq \frac{C_p}\e,
$$
 where $B$ is any closed ball of radius 1, and $C_p$ depends on $p$ but not on $\e$. Hence we have
\begin{equation}
\label{e3.35}
\Vert\nabla u\Vert_\infty\leq\frac{C}\e,
\end{equation}
$C>0$ universal. Hence the condition for $q$ on the road becomes
$$
Dq_{\xi\xi}=\e cq_\xi+\e/\mu q_\zeta=\frac{O(1)}{M_{\e,\lambda}}=o(1).
$$
This implies, in order to keep the positivity of $q$, that $q_\xi(0,0)$ is universally bounded; as a result, there is $\gamma_0>0$ universal, and $\xi_0>0$ universal as well, such that 
$$
\gamma_0\leq q(\xi,0)\leq\frac1{\gamma_0},\quad \xi\in(-\xi_0,\xi_0).
$$
Thus $q(\xi,\zeta\geq\underline q(\xi,\zeta)$ with
$$
\begin{array}{rll}
-d\Delta \underline q+\e c\underline q_\xi+\underline q=&0\quad ((\xi,\zeta)\in(-\xi_0,\xi_0)\times(-1,0))\\
\underline q(\xi,0)=\gamma_0,\ \underline q(\xi,-1)=&0\quad (\xi\in(-\xi_0,\xi_0))\\
\underline q(\pm\xi_0,\zeta)=&0\quad(\zeta\in(-1,0))
\end{array}
$$
and the strong maximum principle implies the existence of a universal $\gamma_0'>0$ such that
$
\underline q(0,\zeta_0)\geq\gamma_0',
$
a contradiction. \hfill$\square$

\noindent Putting Lemmas \ref{l3.2} and \ref{l3.3} together, we obtain the gradient bound in the region where $f_\e(u)$ is active:
\begin{cor}
\label{c3.1}
 There is $M>0$ universal such that $\vert\nabla u(x,y)\vert\leq M$ if $u(x,y)\leq\e$. 
\end{cor}
\noindent{\sc Proof.} Consider $\lambda\in(0,1)$ and a point $(x_0,x_0)\in\overline\Sigma$ such that 
$
u(x_0,y_0)=\lambda\e.
$
By Lemma \ref{l3.1}, it is sufficient to assume $y_0\geq-2\e.$  From Lemmas \ref{l3.2} and \ref{l3.3} we have, in the Lipschitz scaling \eqref{e3.17}:
$$
0\leq p(\xi,\zeta)\leq C,\quad (\xi,\zeta)\in B_2(0)\cap{\zeta<0}.
$$
The Wentzell condition reads once again 
$$
p_{\xi\xi}(\xi,0)=o(1),\ -2\leq\xi\leq2,
$$
so that
$
p_\xi(\xi,0)=o(1),$ for $-2\leq\xi\leq0$.
in order to keep the positivity of $p$. Extend $p(\xi,0)$ in a $C^{1,1}$ fashion inside $[-2,2]\times[-2,0]$, call $\tilde p(\xi,\zeta)$ this extension. Then $p-\tilde p$ solves an elliptic equation with bounded right handisde inside $[-2,2]\times[-2,0]$, while satisfying a Dirichlet condition on the road segment $(-2,2)$. Hence, the classical local elliptic estimates yield
$$
\Vert p\Vert_{L^r(B_1(0)\cap\{\zeta<0\})}\leq C_r,\ \hbox{$C_p>0$ depending only on $r$.}
$$
This implies, in turn, the uniform boundedness of $\vert\nabla p(0,\zeta_0)\vert$, hence the corollary. \hfill$\square$
\subsection{The gradient bound away from the region $u\sim\e$}
Let us mention from the outright that it is the sole place where we will need $D\geq d$. Once again, we do not know if it is essential.

\noindent{\sc Proof Theorem \ref{t4.1}.} We set 
$$
\Omega_\e=\Sigma\backslash\{u\leq\e\}.
$$
Because $\partial_x u<0$, the set $\{u\leq\e\}$ sits on the right of a smooth graph 
$$
\Gamma_\e=\{(h_\e(y),y),\ y\in[-L,0]\},
$$
where $h_\e$ is a smooth function (whose derivatives may nonetheless blow up as $\e\to0$). We may assume $h_\e(0)=0$, and we set $x_\e=h_\e(-L)$; we have thus
$$
\partial\Omega_\e=(-\infty,0)\times\{0\}\cup\Gamma_\e\cup(-\infty,x_\e)\times\{-L\}.
$$
Let us first bound $u_y$ from above. On $\Gamma_\e$ we have $u_y\leq C$. On $(-\infty,x_\e)\times\{-L\}$ we have $u_y=0$. So, let us find a differential inequality for $u_y$ on $(-\infty,0)\times\{0\}$. In $\Omega_\e$ we have
$$
\begin{array}{rll}
-Du_{xx}+cu_x=&D(-u_{xx}+\di\frac{c}du_x)+(1-\di\frac{D}d)cu_x\\
=&Du_{yy}+(1-\di\frac{D}d)cu_x\\
\geq&Du_{yy},
\end{array}
$$
simply because $c>0$, $u_x<0$ and $D\geq d$. This inequality carries over to the boundary $(-\infty,x_\e)\times\{0\}$ to yield
$
D\partial_yu_y+u_y\leq0.$
We also have $-D\Delta u_y+c\partial_xu_y=0$ in $\Omega_\e$, therefore $u_y$ can only assume a positive maximum on $\Gamma_\e$, where it is uniformly bounded from above. If not, it is bounded from above by 0.

\noindent Let us now bound $u_x$ from below. We now know that $u_y\leq C$, so, using $u_x<0$ we write
$$
\begin{array}{rll}
u_{xx}(x,0)-cu_x(x,0)\leq &C\quad(x<0)\\
u_x(0,0)=&O(1)\quad\hbox{because of Cororllary \ref{c3.1}}
\end{array}
$$
Integrating this inequality backward, and using Proposition \ref{p3.3}, we obtain
$
u_x(x,0)\geq-C$, $C>0$ universal.
On $(-\infty,x_\e)\times\{0\}$ we have $\partial_yu_x=0$, and $u_x$ bounded on $\Gamma_\e$. Thus, $u_x$ is bounded from below in $\Omega_\e$.

\noindent Finally, we may bound $u_y$ from below: now that we know the boundedness of $u_x$, the Wentzell condition becomes
$$
D\partial_yu_y(x,0)+u_y(x,0)=O(1).
$$
On $\Gamma_\e$, $u_y$ is still bounded, while it still satisfies the Dirichlet condition on $(-\infty,x_\e)$. The lower bound follows from the maximum principle.
\hfill$\square$

\section{The free boundary near $(0,0)$ in the single species model}
  For the reader's convenience, we recall the system solved by $(c,\Gamma,u)$:
\begin{equation}
\label{e4.3}
\left\{
\begin{array}{rll}
-d\Delta u+c\partial_x u=&0\quad(x,y)\in\{u>0\}\\
u_\nu=&0\quad((x,y)\in\Gamma\\
\ \\
-Du_{xx}+c\partial_xu+1/\mu u_y=&0\quad\hbox{for $x\in\RR$, $y=0$}\\
u_y(x,-L)=&0\\
\ \\
u(-\infty,y)=&1,\quad u(+\infty,y)=0.
\end{array}
\right.
\end{equation} 
Let us first explain the heuristics of Theorem \ref{t2.4}.  The starting point is a very simple global asymptotic solution, in the limit $\e\to0$ of the parabolic scaling: 
\begin{equation}
\label{e2.4}
u(x,y)=\e v(\frac{ x}{\sqrt\e},\frac{y}\e),
\end{equation}
in other words, the following situation:

\includegraphics{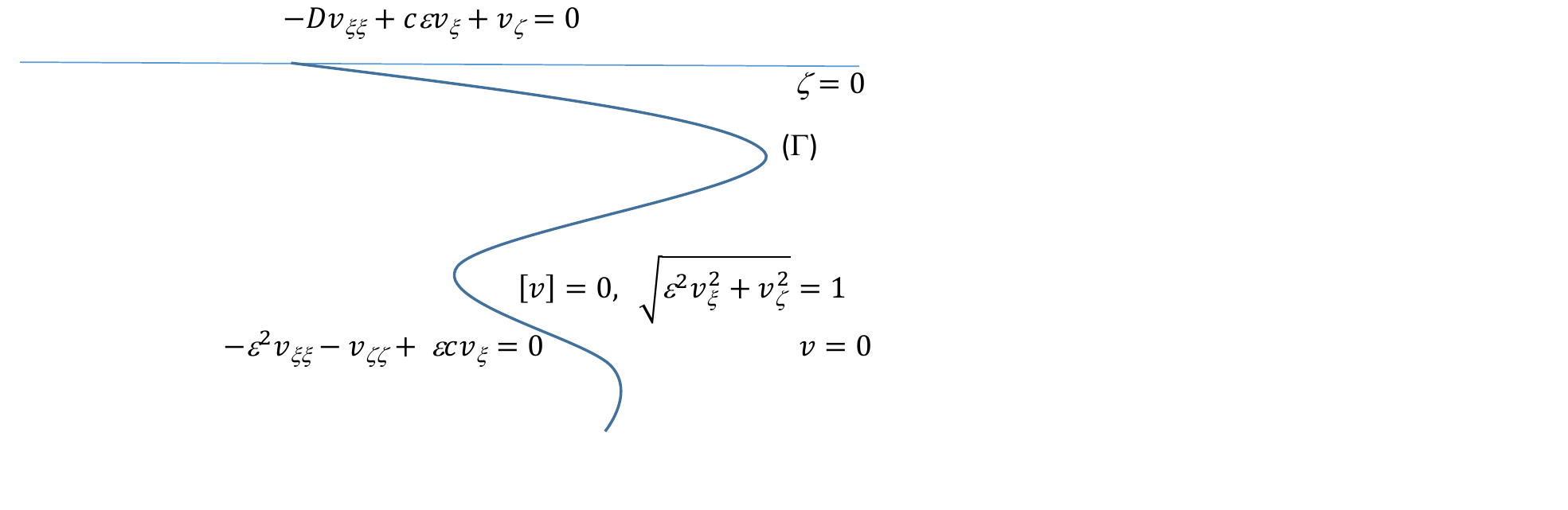}
Setting $\e=0$ yields $v_{\zeta\zeta}\equiv0$, hence $v(\xi,0)\equiv1$ for $\xi<0$. So, we have $v(\xi,0)=\di\frac{\xi^2}{2D}$ and
$$
v(\xi,\zeta)=\biggl(\zeta+\frac{\xi^2}{2D}\biggl)_+,
$$
hence the free boundary has the form \eqref{e2.1000}. So, we are going to make this rigorous in the remaining of the section.

\noindent We may assume that $\Gamma$ meets the axis $\{y=0\}$ at the point (0,0), and start the proof of Theorem \ref{t2.2}. If $\bar X=(\bar x,\bar y)\in\Sigma$ is a free boundary point, and  $\delta>0$ is small enough so that $B_{2\delta}(\bar X)\subset\Sigma$, rescale $u$ at scale $\delta$:
\begin{equation}
\label{e4.1}
x-\bar x=\delta\xi,\ y-\bar y=\delta\zeta,\quad u(x,y)=\delta u_\delta(\xi,\zeta).
\end{equation}
Let $\Gamma_\delta$ be the rescaled free boundary. As is familiar we have
\begin{equation}
\label{e4.2}
\begin{array}{rll}
-d\Delta u_\delta+\delta c\partial_\xi u_\delta=&0\quad (B_2(0))\\
\vert\nabla u_\delta\vert=&1\quad(\Gamma_\delta)\\
u(0)=&0.
\end{array}
\end{equation}
The uniform gradient bound comes from Theorem \ref{t2.1}.
We may therefore safely forget about $u_\e$, and concentrate on the solution $u(x,y)$ of \eqref{e1.2}. And so, in the sequel, we will use the letter $\e$ for any small parameter, without any further reference to the approximate solution constructed in the preceding section.

\noindent This section more or less follows the organisation of Section 3: first, we will state a (more standard) rigidity property adapted to the needs of the asymptotic situations that we will encounter. We will then prove theorem \ref{t2.4}.
\subsection{Another 2D rigidity result}
The following proposition is a consequence of the 2D monotonicity formula.
\begin{prop}
\label{p5.100}
Consider a  solution $u(x,y)$ of
\begin{equation}
\label{e5.601}
\begin{array}{rll}
\Delta u=&0, \quad (\{u>0\})\\
\vert\nabla u\vert=&1\quad (\partial\{u>0\}),\\
u(x,0)=&0\quad (x\in\RR).\\
\end{array}
\end{equation} 
Assume that $\partial\{u>0\}$ meets the $x$-axis at the point (0,0), which s the only point where it may not be analytic. Then $u(x,y)=y^-$.
\end{prop} 
\noindent{\bf Proof.} It is based on \cite{CafS}, Theorem 12.1. Set 
$$
u_1(x,y)=u(x,y),\quad u_2(x,y)=u(x,-y).
$$
The functions $u_1$ and $u_2$ are harmonic in disjoint domains, and have a common zero line: the axis $\{y=0\}$. Thus (\cite{CafS}, Theorem 12.1) the quantity
$$
\frac1{R^4}\int_{B_R(0)}\vert\nabla u_1\vert^2dxdy\int_{B_R(0)}\vert\nabla u_2\vert^2dxdy
$$
is increasing in $R$. As a consequence, the quantity
$$
J_u(R)=\frac1{R^2}\int_{B_R(0)}\vert\nabla u\vert^2dxdy
$$
is also increasing in $R$. Because of the gradient bound and nondegeneracy, the quantities
$
\di\lim_{R\to0}J_u(R)$ and $\di\lim_{R\to+\infty}J_u(R)
$
exist and are nonzero. Notice that $J_u$ is invariant under Lipschitz scaling, that is
$$
J_u(R)=\int_{B_1(0)}\vert\nabla u_R\vert^2dxdy=J_{u_R}(1),\quad u_R(x,y)=\frac1R u(Rx,Ry).
$$ 
The analyticity of the free boundary entails enough compactness to infer that the family $(u_R)_R$ converges in $H^1(B_1(0))$, as $R\to0$ (resp. $R\to+\infty$), possibly up to a subsequence, to a limit $u^-$ (resp. $u^+$). Another application of Theorem 12.1 of \cite{CafS} to $u^\pm$ yields
$$
J_{u^\pm}(r)=\mathrm{Constant},\ \hbox{hence}\ J'_{u^\pm}(r)=0.
$$
Hence we have
$$
\lim_{R\to0} J_u'(R)=0,\quad \lim_{R\to+\infty}J_u'(R)=0.
$$
Following the proof of the theorem, we have that $\lim_{R\to0} J_u'(R)$ (resp. $\lim_{R\to+\infty} J_u'(R)$)  is proportional to
$$
\di{\frac{\di\int_{\partial B_1(0)}\vert\nabla u^-\vert^2dxdy}{\di\int_{B_1(0)}\vert\nabla u^-\vert^2dxdy}}-2\quad\biggl(\hbox{resp.} \di\frac{\di\int_{\partial B_1(0)}\vert\nabla u^+\vert^2dxdy}{\di\int_{B_1(0)}\vert\nabla u^+\vert^2dxdy}-2\biggl).
$$
However, still following the proof of the theorem, we find out that the above quantities are zero if and only if
$$
u^\pm(y)=y^-;
$$
actually $u^\pm$ have to be proportional to $y^-$, but the free boundary relation imposes the proportionality coefficient. Thus $J$ has the same limits at 0 and $+\infty$, thus $J$ is constant. Thus $J'(R)\equiv0$, and 
$
u(x,y)=y^-.
$ 
This proves the proposition. \hfill$\square$
\subsection{Analysis of the free boundary}
The first step is to show that the Lipschitz scale inside breaks down as the free boundary approaches the horizontal axis, and becomes flatter and flatter.
\begin{lm}
\label{l4.1}
There is $C>0$ universal such that, for $\e>0$ small enough, we have
$$
\Gamma\cap B_\e(0)\subset (-\e,0]\times(-C\e^2,0].
$$
\end{lm}
\noindent{\sc Proof.} Assume the lemma to be false, that is, there is a sequence $(x_\e,\delta_\e)$ with
$$
\Gamma\cap\partial B_\e(0)=\{(x_\e,\delta_\e)\},
$$
with
\begin{equation}
\label{e4.4}
\lim_{\e\to0}x_\e\sqrt{\delta_\e}=0.
\end{equation}
Drop the subscript ${}_\e$ and scale with $\delta$ as in \eqref{e4.1}, with $\bar X=(0,0)$. The Wentzell condition implies that 
$$
u_\delta(\xi,0)=O(\delta\xi^2)\to_{\delta\to0}0.
$$
The situation implies the existence of $\bar X_\delta$ in $B_1(0)\cap\{u_\delta>0\}$ whose distance to $\Gamma_\delta$, as well as to $\{y=0,\}$, is universally
controlled from below. From nondegeneracy, $u_\delta(\bar X_\delta)$ is universally controlled from below. We send $\delta$ to 0 and use the compactness of $(u_\delta)_\delta$ provided by the uniform gradient bound: we recover a function $u_\infty(\xi,\zeta)$, as well as an asymptotic smooth (possibly outside the origin) free boundary $\Gamma_\infty$. Call $\Omega_\infty$ the area limited by $\{\zeta=0\}$ and $\Gamma_\infty$, we have
\begin{equation}
\label{e4.5}
-\Delta u_\infty=0\ \hbox{in $\Omega_\infty$},\quad u_\infty=0\ \hbox{on $\partial\Omega_\infty$.}
\end{equation}
 nondegeneracy implies that $u_\infty(\xi,-1)$ is uniformly controlled from below, which implies in turn
$$
u_\zeta(\xi,0)\leq -2q
$$
for some universal $q$.
Returning to $u_\delta$ we obtain the same kind of bound: for every $\rho>0$, there is $q_\rho>0$ such that
$$
\partial_\zeta u_\delta(\xi,\zeta)\leq -q_\rho\ \hbox{for small $\delta$, $-\xi_\delta<\xi\leq-\rho$.}
$$
We have, because of \eqref{e4.4}:
$$
\lim_{\delta\to0}\xi_\delta=+\infty.
$$
The Wentzell condition at $\{\zeta=0\}$ implies
$$
-D\partial_{\xi\xi}u_\delta+\delta c\partial_\xi u_\delta\geq q\delta(\un_{(-\xi_\delta,-\rho)}-C\un_{(-\rho,0)}),\ u_\delta(0,0)=\partial_\xi u_\delta(0,0)=0.
$$
Integration of this inequality yields (this is by now a routine argument) yields
$$
u_\delta(-\xi_\delta,0)\leq -\frac{q_\rho}3\delta<0,
$$
the sought for contradiction. This proves the lemma. \hfill$\square$

The next lemma shows that, if the free boundary wiggles before reaching the point (0,0), it does so in a controlled fashion.
\begin{lm}
\label{l4.2}
For every $\e>0$, let $\zeta_\e^\pm$ be defined as
$$
\zeta_\e^-=\inf\{\zeta\in\RR_+: (-\e,\zeta)\in\Gamma\},\quad \zeta_\e^+=\sup\{\zeta\in\RR_+: (-\e,\zeta)\in\Gamma\} 
$$
There is $q\in(0,1)$, universal, such that 
$$
\vert\zeta_\e^-\vert\geq q\vert\zeta_\e^+\vert.
$$
\end{lm}
\noindent{\sc Proof.} 
Do the Lipschitz scaling \eqref{e4.1} with $\bar X=(-\e,0)$ and $\delta=\zeta_\e^+$. Assume that, for a sequence $\delta_n\to0$ we have:
$$
\lim_{n\to+\infty}\inf\{\zeta>0:\ (0,-\zeta)\in\Gamma_{\delta_n}\}=0.
$$
This implies immediately
$$
\lim_{n\to+\infty}u_{\delta_n}(\xi,0)=0,\quad\xi\geq0,\quad
\lim_{n\to+\infty}\partial_\xi u_{\delta_n}(0,0)=0.
$$
By compactness, we obtain a nontrivial limiting couple $(\Gamma_\infty,u_\infty)$ which solves the free boundary  problem in $\RR\times\RR_-$, and such that $(0,0)\in\Gamma_\infty$. For $\xi<0$ we have $\partial_{\xi\xi}u_\infty(\xi,0)=0$, therefore $u_\infty(\xi,0)\equiv0$ for $\xi\in\RR$. This is against Proposition 
\ref{p5.100}. 
\hfill$\square$

\noindent Lemma \ref{l4.2} will trigger an obvious analogy between the equation for $y$ and the obstacle problem, in the sense that $u$ behaves in a quadratic fashion in the vicinity of the point where the inside free boundary hits the axis $\{y=0\}$. This is expressed in the next proposition.

\begin{prop}
\label{p4.1}
For some$C>0$,  we have, for $x\in[-1,0]$:
$
\di\frac{x^2}C\leq u(x,0)\leq Cx^2.
$
\end{prop}
\noindent{\sc Proof.} Of course, it is enough to prove the result for small $x$. Consider any sequence $(\e_n)_n$ going to 0,
without loss of generality we may shift the origin from $-\e_n$ to 0.  
Let us set 
$$\delta_n:=\zeta_{\e_n}^+.
$$
We wish to prove that $u(-\e_n,0)$ is of the order $\delta_n$; because of the gradient bound it is certainly an $O(\delta_n)$. Let us prove the converse, i.e. that
$$
\delta_n=O\biggl(u(-\e_n,0)\biggl).
$$
Assume this is not true.
Do, as in the preceding situation, the Lipschitz scaling \eqref{e4.1} with $\delta=\delta_n$. Let $(\Gamma_n,u_n)$ be the scaled versions of $u$ and $\Gamma$. Send $n\to+\infty$,  we recover an unbounded analytic curve $\Gamma_\infty=\{h_\infty(\zeta),\zeta\}$, lying strictly below $\{\zeta=0\}$, and a positive harmonic function $u_\infty$, defined above $\Gamma_\infty$, satisfying the free boundary relation on $\Gamma$. This is impossible, due once again to Proposition \ref{p5.100}. The same proposition will imply 
$$
\partial_\zeta u_\infty(0,0)=1,
$$
so that, for the unscaled function $u(x,y)$, there is $q>0$, universal, such that $u_y(x,0)\geq q$ for $x<\in(-\e_n,0)$ as soon as $n$ is large enough, except perhaps in an $o(\e_n)$-neighbourhood of 0. So, by scaling, we have $u_y(x,0)\geq q$ as $x<0$ is small enough. Integration of the Wentzell boundary condition yields the quadratic nondegeneracy of $u$ on the road, in a small neighbourhood of 0.
\hfill$\square$

\noindent The last ingredient that we need is that, at each point of the free boundary near $(0,0)$, the normal is almost vertical.
\begin{lm}
\label{l4.6}
There is $\delta_0>0$ such that, if $x\in(-\delta_0,0)$, there is a unique $y=k(x)$ such that $(x,y)\in\Gamma$. For $y\in(-\delta_0,0)$, let $\nu(x)=(\nu_1(x),\nu_2(x))$ be the outward normal to $\Omega$ at the point $(x,k(x))$. Then we have, for some $C>0$:
\begin{equation}
\label{e4.72}
-Cx^2\leq k(x)\leq -\frac{x^2}C
\end{equation}
Moreover we have:
$$
\lim_{x\to0}\nu_1(x)=0,\quad\lim_{x\to0}\nu_2(x)=-1.
$$
\end{lm}
\noindent{\sc Proof.} For every $\e$ small, let $\delta_\e$ be the smallest $\delta$ such that 
$$
(-\e,-\delta)\in\Gamma,
$$
from Lemma \ref{l4.2} the largest $\delta$ such that this property holds is also of the order $\delta_\e$. We do the Lipschitz scaling with 
$$
\bar X=(-\e,0),\quad\delta=\delta_\e.
$$
Let $u_\e$ be the rescaled function, and send $\e$ to 0. From the above lemmas we obtain a couple $(\Gamma_\infty,u_\infty)$ where $\Gamma_\infty$ is a union of analytic curves, one above the other,  trapped in a bounded strip of the form $\RR\times(-M,0)$, and $u_\infty$ solves the free boundary problem in the set $\Omega_\infty$ between $\Gamma_\infty$ and $\{\zeta=0\}$. Moreover, the top part of $\partial\Omega_\infty$ is a straight line, and there is $q>0$ such that 
$$
u_\infty(\xi,0)=q.
$$
We also know that $\partial_\xi u_\infty\leq0$, therefore it has two limits $u_\infty^\pm(\zeta)$ as $\zeta\to\pm\infty$. Hence these limits are non trivial, and they also solve the free boundary problem. So, an easy computation yield
$$
u^{\pm}_\infty(\zeta)=(q+\zeta)^+,\quad\zeta<0.
$$
Thus we have
$$
u_\infty(\xi,\zeta)=(q+\zeta)^+,\quad\Gamma_\infty=\{\zeta=-q\}.
$$
Thus, the whole family $(u_\e)_\e$ converges to $(q+\zeta)^+$. This also implies that, for $\xi\in(-1,1)$, the uniqueness of $\zeta$ such that $(\xi,\zeta)\in\Gamma$. Hence $\Gamma_\e$ is an analytic graph $(\xi,k_\e(\xi))$ with $(k_\e')_\e$ bounded from Lemma \ref{l4.2}. Hence we have
$$
\lim_{\e\to0}k_\e(x)=-q,\ \lim_{\e\to0}k_\e'(x)=0\quad\hbox{uniformly in $\xi\in(-1,1)$.}
$$
This implies the convergence of the normals, as stated in the lemma. \hfill$\square$

\noindent {\sc Proof of Theorem \ref{t2.2}.}
We finally revert to the Lipschitz-parabolic scaling, and set:
\begin{equation}
\label{e4.70}
x=\sqrt\e \xi,\ y=\e \zeta,\quad u_\e(\xi,\zeta)=\frac1\e u(\sqrt\e\xi,\e\zeta).
\end{equation}
Let $\Gamma_\e$ be, once again, the free boundary of $u_\e$. It is enough to prove that, for small $\e>0$, the only $\zeta<0$ such that 
$
(-1,\zeta)\in\Gamma_\e,
$ 
then we have 
$$
\zeta=-\frac1{2D}.
$$
Pick any small $\delta\in(0,2)$. From Lemma \ref{l4.6}, for $-2\leq\xi\leq-\delta$, there is a unique $\zeta_\e(\xi)\leq0$ such that
$$
(\xi,\zeta_\e(\xi))\in\Gamma_\e.
$$
From Lemma \ref{l4.6} we have 
$$
\lim_{\e\to0}\partial_\zeta u_\e(\xi,\zeta_\e(\xi))=\lim_{\e\to0}\partial_yu(\sqrt\e\xi,\e\zeta_\e(\xi))=1.
$$
From the proof of the same lemma we have
$$
\lim_{\e\to0}\partial_\zeta u_\e(\xi,0)=1.
$$
The Wentzell condition for $u_\e$ is 
$$
\begin{array}{rll}
D\partial_{\xi\xi}u_\e-\e c\partial_\xi u_\e=&\partial_\zeta u_\e\sim_{\e\to0}1\quad\hbox{for $\xi\leq-\delta$},\\
u_\e(0,0)=&\partial_\xi u_\e(0,0)=0.
\end{array}
$$
Integrating this relation and taking the boundedness of $\partial_\zeta u_\e$ into account yields
$$
\lim_{\e\to0}u_\e(-1,0)=-\frac1{2D}.
$$
And, from the proof of Lemma \ref{l4.6}, we have
$$
\lim_{\e\to0}u_\e(-1,\zeta)=(\zeta+\frac1{2D})^+.
$$
This implies the result. \hfill$\square$

\begin{rem} 
\label{r5.1}
The gradient bounds, as well as the upper and lower bounds on $c$, do not depend on $L$. In the proof of Theorem \ref{t2.2}, the convergence of  $u_\e$ to $(\zeta+\di\frac1{2D})^+$ is uniform in $L\geq 1$.
\end{rem}

\noindent {\bf Acknowledgement.} L.A. Caffarelli is supported by NSF grant DMS-1160802. The research of J.-M. Roquejoffre  has received funding from the ERC under the European Union’s Seventh Frame work Programme (FP/2007-2013) / ERC Grant Agreement 321186 - ReaDi. He also acknowledges  J.T. Oden fellowships, for visits at the University of Texas. 

\noindent 
{\footnotesize
 
}


\begin{thebibliography}{99}
\bibitem{ACF} {\sc H.W. Alt,  L. A. Caffarelli,} {\it Existence and regularity for a minimum problem with free boundary,} Journal für die reine und angewandte Mathematik {\bf 325} (1981), 105--144.
\bibitem{BCN}\textsc{H. Berestycki, L.A. Caffarelli, L. Nirenberg}, {\it Uniform estimates for the regularization of free boundary problems,} Analysis and Partial Differential Equations, 567--619, Lecture Notes in Pure and Applied Mathematics, {\bf 122}, Dekker, New York, 1990.
\bibitem{BNS} {\sc H. Berestycki, B. Nicoalenko, B. Scheurer,} {\it . Traveling wave solutions to combustion models and
their singular limits}, SIAM J. Math. Anal. {\bf 16} (1985), 1207--1242.
\bibitem{BRR2}
{\sc H.~Berestycki, J.-M. Roquejoffre,  L.~Rossi},
{\em The influence of a line with fast diffusion on {F}isher-{KPP}
  propagation.} J. Math. Biol. {\bf 66} No. 4-5 (2013), 743--766.
\bibitem{BRR4} {\sc H. Berestycki, J.-M. Roquejoffre, L. Rossi} {\it The shape of expansion induced by a line with fast diffusion in Fisher-KPP equations}, Communication in Mathematical Physics, {\bf 343} (2016), 207--232.
\bibitem{CafS}{\em } \textsc{L.A. Caffarelli, S. Salsa}, {\it A geometric approach to free boundary problems, } Graduate Studies in Mathematics, 68, American Math. Soc., Providence, R.I., 2005.
\bibitem{ACC} {\sc A.-C. Coulon Chalmin}, {\it Fast propagation in reaction-diffusion equations with fractional diffusion}, PhD thesis, 2014. Online manuscript: thesesups.ups-tlse.fr/2427/
\bibitem{Di1} {\sc L. Dietrich}, {\it Existence of Travelling Waves for a Reaction–Diffusion System with a Line of Fast Diffusion}, Appl. Math. Research Express, {\bf 2} (2016), 204--252.
\bibitem{DiR} {\sc L. Dietrich, J.-M. Roquejoffre}, {\it Front propagation directed by a line of fast diffusion: large diffusion and large time asymptotics
}, J. Ecole Polytechnique., {\bf 4} (2016), 141--176.
\bibitem{LTru} {\sc Y. Luo, N. S. Trudinger,}
{\it Linear second order elliptic equations with Venttsel boundary conditions,} Proc. Roy. Soc. Edinburg Sect. A, {\bf 118}
(1991),  193--270.
\bibitem{ZBLM}{\sc Ya. B. Zel'dovich, G. Barenblatt G., V.B. Librovich, G.M. Makhviladze}. {\it The Mathematical Theory of Combustion and Explosions}, Consultants Bureau, 1985.
\end{thebibliography}
\end{document}